\documentclass{ifacconf}

\usepackage{etoolbox}

\makeatletter
\pretocmd{\@ssect}{\NR@gettitle{#6}}{}{}
\makeatother

\usepackage{subcaption}

\makeatletter
\@ifundefined{c@part}{\newcounter{part}}{}
\@ifundefined{thepart}{\newcommand\thepart{0}}{}
\@ifundefined{theHpart}{}{}
\let\old@ssect\@ssect % Store how ifacconf defines \@ssect
\makeatother

\usepackage{natbib}
\usepackage{hyperref}
\hypersetup{
    colorlinks = true,
    linkcolor = {blue},
    citecolor = {blue}
}

\makeatletter
\def\@ssect#1#2#3#4#5#6{%
  \NR@gettitle{#6}% Insert key \nameref title grab
  \old@ssect{#1}{#2}{#3}{#4}{#5}{#6}% Restore ifacconf's \@ssect
}
\makeatother

\usepackage{amsmath, amsfonts, amssymb, bm}
\usepackage{cleveref}
\usepackage[dvipsnames]{xcolor}
\usepackage{graphicx} 

\newcommand{\cA}{\mathcal{A}}

\newcommand{\RR}{\mathbb{R}}

\usepackage[ruled,vlined,algo2e]{algorithm2e}
\begin{document}
\begin{frontmatter}

\title{Scalable Method for Mean Field Control with Kernel Interactions via\\ Random Fourier Features}%\thanksref{footnoteinfo}} 

% \thanks[footnoteinfo]{Sponsor and financial support acknowledgment goes here. Paper titles should be written in uppercase and lowercase letters, not all uppercase.}

\author[First]{Zhongyuan Cao} 
\author[Second]{Kaustav Das} 
\author[Third]{Nicolas Langren\'e}
\author[Fourth]{Mathieu Lauri\`ere}

\address[First]{NYU-ECNU Institute of Mathematical Sciences, NYU Shanghai, Shanghai, 200126, People’s Republic of China (e-mail: zc3151@nyu.edu).}
\address[Second]{School of Mathematics, Monash University, Clayton, 3800, Australia (e-mail: kaustav.das@monash.edu)}
\address[Third]{Guangdong Provincial/Zhuhai Key Laboratory of Interdisciplinary Research and Application for Data Science, Department of Mathematical Sciences, Beijing Normal-Hong Kong Baptist University, Zhuhai, 519087, People’s Republic of China (e-mail: nicolaslangrene@bnbu.edu.cn)}
\address[Fourth]{Shanghai Center for Data Science; NYU-ECNU Institute of Mathematical Sciences, NYU Shanghai, Shanghai, 200126, People’s Republic of China (e-mail: ml5197@nyu.edu).}

% \begin{abstract}                % Abstract of 50--100 words
% \color{red}
% These instructions give you guidelines for preparing papers for IFAC
% technical meetings. Please use this document as a template to prepare
% your manuscript. For submission guidelines, follow instructions on
% paper submission system as well as the event website.
% \end{abstract}

\begin{abstract}
We develop a scalable algorithm for mean field control problems with kernel interactions by combining particle system simulations with random Fourier feature approximations. The method replaces the quadratic-cost kernel evaluations by linear-time estimates, enabling efficient stochastic gradient descent for training feedback controls in large populations. We provide theoretical complexity bounds and demonstrate through crowd motion and flocking examples that the approach preserves control performance while substantially reducing computational cost. The results indicate that random feature approximations offer an effective and practical tool for high dimensional and large scale mean field control.
\end{abstract}


\end{frontmatter}
%===============================================================================

\section{Introduction}

Many systems involving large populations of interacting agents can be described by stochastic dynamics coupled through nonlocal terms. Examples include models of pedestrian motion, coordinated motion in animal groups, and more broadly multi-agent systems in which the behavior of each individual depends on an aggregate of the surrounding population. These interactions are often spatial or similarity based and can be naturally represented through kernel functions. When such systems must be controlled, one is led to mean field control problems with nonlocal interactions.

Mean field control (MFC) provides a mathematical framework for optimizing the collective behavior of a continuum of agents; it differs from the mean field game setting where the agents are non-cooperative and play a Nash equilibrium. We refer e.g. to the monographs~\citep{bensoussan2013mean} and~\citep{carmona2018probabilistic} for more details. It can be interpreted both as a stand-alone stochastic control problem in which the law of the state enters the dynamics and cost, and as the limit of large $N$-agent optimal control problems. In this limit, the optimal control becomes decentralized: each agent uses the same feedback rule depending only on its own state. As far as non-local interactions are concerned, most works have focused on interactions through the first moment in the form of linear-quadratic models, see e.g.~\citep{caines2014mean,ni2014indefinite,bensoussan2016linear,graber2016linear,yong2017linear}. When interactions occur through kernels, the controlled drift and cost depend on the convolution of the population distribution with a prescribed kernel that encodes the influence of nearby agents.

Numerical methods for MFC often rely on particle approximations of the population distribution~\citep{fouque2020deep,carmona2022convergence,dayanikli2024deep}. In this setting, evaluating a nonlocal interaction term such as $(K \star \mu)(x)$ requires computing all pairwise kernel evaluations $K(x_i - x_j)$ for $i,j$ in $[N]$, where $N$ is the number of particles. This leads to a computational cost of order $N^2$ per time step, which becomes prohibitive when $N$ is large or when long time horizons are considered. Grid-based or PDE-based approaches face other difficulties, including dimensionality and the cost of discretizing nonlocal operators. Reducing the computational effort associated with kernel interactions is therefore an important step toward scalable algorithms for MFC.

Random Fourier Features (RFF), introduced by \citet{rahimi2007random},  provide an efficient way to approximate translation-invariant positive definite kernels. By sampling from the spectral density of the kernel, one constructs a finite-dimensional random feature map $\Phi$ such that $K(x-y)$ is approximated by $\Phi(x)^{\top}\Phi(y)$. Replacing the kernel by this approximation transforms the convolution into an inner product between $\Phi(x_i)$ and an aggregated feature vector. This reduces the computational complexity of evaluating the interaction term from order $N^2$ to order $NM$, where $M$ is the number of random features and can be taken much smaller than $N$. The resulting approximation error decreases as $M^{-1/2}$ and the method applies to many kernels used in mean field models, including Gaussian, Matérn, and Cauchy kernels, see~\citet{langrene2025mixture}.

The main purpose of this paper is to leverage RFF to obtain a scalable numerical method for MFC problems with nonlocal interactions. We show how to incorporate the RFF approximation directly into particle-based simulations of MFC dynamics. This provides efficient approximations of kernel convolutions and related quantities that appear in controlled systems. We combine this idea with a deep learning method for MFC proposed in~\citep[Algorithm 1]{carmona2022convergence} which is able to handle high-dimensional problems. We illustrate the method on two representative models: a pedestrian model with Gaussian or Matérn interactions and a flocking model with Cauchy interactions. In both cases, the RFF-based approximation leads to significant computational gains while retaining sufficient accuracy for control and simulation.

Previously, only a few works used kernel methods combined with MFC. \cite{liu2021computational} studied an approach relying on kernel-based representations of mean field interactions and feature-space expansions in the spirit of kernel methods in machine learning but only for Gaussian kernels.  \cite{mou2022numerical} approximates the solution of a mean-field game using a Gaussian process, the inference of which is accelerated by the random Fourier features approximation of the covariance kernel. \cite{agrawal2022random} proposed an efficient solution approach for high-dimensional nonlocal mean-field game systems by approximating their Gaussian kernel interaction function by random features.

The paper is organized as follows. Section~\ref{sec:MFC} introduces MFC problems with kernel interactions, presents motivating examples, and describes the $N$-particle approximation. Section~\ref{sec:RFF} explains the RFF methodology and its application to the convolution terms arising in MFC. Section~\ref{sec:appli} describes how to apply our method to two examples. Section~\ref{sec:numerics} reports numerical experiments for the two models. Concluding remarks are given in Section~\ref{sec:ccl}.

\section{Mean Field Control Problems with Kernel Interactions}
\label{sec:MFC}
In this section, we introduce a generic MFC problem and specialize it to kernel interactions.

\subsection{Mean Field Control Problem}
Let $T\in[0,\infty)$ be a fixed time horizon.  Let $\cA$ be an appropriate set of controls that we will consider for the mean field control problem. The action set is denoted by $\mathbb{A}$, which usually is a subset of $\mathbb{R}^k$. A control is a $\mathbb{A}$-valued stochastic process $\bm{\alpha} = (\alpha_t)_{t \in [0,T]}$ on $[0,T]$. In the mean field control problem, we seek an optimal control for a population of agents, where each agent follows the same control rule. For a given control ${\bm{\alpha}}\in\cA$, the individual state of a representative agent has the following dynamics evolving in continuous time:
\begin{equation}\label{eq:1}
dX^{\bm{\alpha}}_t = b(t, X^{\bm{\alpha}}_t, \alpha_t, \mu^{\bm{\alpha}}_t)\,dt + \sigma\,dW_t,
\end{equation}
where $\bm{W}=(W_t)_{t \in [0,T]}$ is a standard $m$-dimensional Brownian motion,
% \item $\alpha_t \in A \subset \mathbb{R}^k$ is the action at time $t$,
$X^{\bm{\alpha}}_t \in \mathbb{R}^d$ is the state at time $t$, under control ${\bm{\alpha}}$,
$\mu^{\bm{\alpha}}_t = \mathcal{L}(X^{\bm{\alpha}}_t)$ is the law of $X^{\bm{\alpha}}_t$, which represents the state distribution of the whole population at time $t$, under control ${\bm{\alpha}}$,
$\sigma$ is a constant $d\times m$ volatility matrix, and $b:[0,T]\times \mathbb{R}^d\times \mathbb{A}\times\mathcal{P}_2(\RR^d)\to \RR^d $ is the drift, where $\mathcal{P}_2(\RR^d)$ denotes the set of probability measures with finite second moment on $\RR^d$. 
% \end{itemize}

Hence, a control $\bm{\alpha}$ induces a probability measure flow 
$\bm{\mu}^{\bm{\alpha}} = (\mu_t^{\bm{\alpha}})_{t \in [0,T]}$, where $X_t^{\bm{\alpha}}$ is the solution of the controlled stochastic differential equation (SDE)~\eqref{eq:1}.

The objective functional that we want to minimize is
\begin{equation}\label{eq:2}
J(\bm{\alpha}) = 
\mathbb{E}\left[\int_0^T f(t, X_t^{\bm{\alpha}}, \alpha_t, \mu_t^{\bm{\alpha}})\,dt
+ g(X_T^{\bm{\alpha}}, \mu_T^{\bm{\alpha}})\right],
\end{equation}
where $f:[0,T]\times \mathbb{R}^d\times \mathbb{A}\times\mathcal{P}_2(\RR^d)\to \RR $ is the running cost, $g:\RR^d\times \mathcal{P}_2(\RR^d)\to \RR$ the terminal cost, and the expectation is taken with respect to the Brownian motion. 

The \textbf{Mean Field Control (MFC) problem} is to seek $\bm{\alpha}^\ast\in\cA$ such that
\begin{equation}\label{eq:3}
J(\bm{\alpha}^\ast)=\inf_{\bm{\alpha}\in \cA} J(\bm{\alpha}),
\end{equation}
that is, we seek an optimal control $\bm{\alpha}^\ast\in \cA$, which can be a feedback control or open-loop control, depending on the admissible control set $\cA$, such that under control $\bm{\alpha}^\ast$, the total expected cost for the population attains the minimum.  
Note that, in contrast with Mean Field Games (MFG), the mean field $\bm{\mu}^{\bm{\alpha}}$ is directly determined by $\bm{\alpha}$—there is no fixed-point consistency condition.

\subsection{Interactions Through a Kernel}
In MFC problems, the interactions between agents are encoded in the population distribution, which can capture the interaction preference when agents' states locate differently. This is often characterized by kernel functions, which we therefore call \textbf{kernel interactions}.  

Among agents, the interactions are mediated by a kernel $K:\mathbb{R}^d\to\mathbb{R}$.
The drift, running cost, and terminal cost depend on the mean field only through its convolution with $K$, namely:
% \begin{align}
% b(t,x,a,\mu) &= \tilde{b}\bigl(t,x,a,(K\star\mu)(x)\bigr), \\
% f(t,x,a,\mu) &= \tilde{f}\bigl(t,x,a,(K\star\mu)(x)\bigr), \\
% g(x,\mu) &= \tilde{g}\bigl(x,(K\star\mu)(x)\bigr),
% \end{align}
$b(t,x,a,\mu) = \tilde{b}\bigl(t,x,a,(K\star\mu)(x)\bigr),$
$f(t,x,a,\mu) = \tilde{f}\bigl(t,x,a,(K\star\mu)(x)\bigr),$
$g(x,\mu) = \tilde{g}\bigl(x,(K\star\mu)(x)\bigr),$ 
for some functions $\tilde{b}, \tilde{f}, \tilde{g}$, where the convolution is
\[
(K\star\mu)(x) = \int_{\mathbb{R}^d} K(x-y)\,\mu(dy).
\]
This convolution represents the local aggregate influence of the population around position $x$.

\subsection{Example 1: Pedestrian Model}
\label{sec:model-crowd}

We consider a model describing the movements of a crowd of pedestrians. We have $d = 2$ and $m = 2$. The state $x\in\mathbb{R}^2$ represents the pedestrian’s position. The dynamics are governed by the stochastic differential equation $dX_t = b(t,X_t,\alpha_t,\mu_t)\,dt + \sigma\,dW_t$, where $b(t,x,a,\mu) = a$. The running cost is given by $f(t,x,a,\mu) = \frac{1}{2}\|a\|^2 + \phi\bigl((K\star\mu)(x)\bigr)$, where $\phi$ is an increasing function penalizing crowded areas. The terminal cost is $g(x,\mu) = \|x - x_{\mathrm{target}}\|^2$, favoring arrival near a preferred destination. The interaction kernel is defined as the generalized Matérn kernel  
\begin{equation}\label{eq:matern}K(x)=\frac{\left(\sqrt{2\beta_2}\left(\frac{\left\Vert x\right\Vert }{\sigma_{K}}\right)^{\frac{\beta_1}{2}}\right)^{\beta_2}}{\Gamma(\beta_2)2^{\beta_2-1}}\mathcal{K}_{\beta_2}\!\left(\sqrt{2\beta_2}\left(\frac{\left\Vert x\right\Vert }{\sigma_{K}}\right)^{\frac{\beta_1}{2}}\right)\end{equation}
with $\beta_1\in(0,2]$ and $\beta_2>0$ \citep{langrene2025mixture}, where $\mathcal{K}_{\beta_2}$ is the modified Bessel function \citep[10.25]{dlmf}. This family of kernels include the Gaussian kernel $K(x) = \exp\!\left(-\frac{\|x\|^2}{2\sigma_K^2}\right)$ as the limit case $\beta_1=2$, $\beta_2\rightarrow\infty$, in which case each pedestrian interacts mostly with nearby agents within range $\sigma_K$. Similar models have been considered, e.g., in~\citep{achdou2015system}.

\subsection{Example 2: Flocking Model}
\label{sec:model-flocking}

We then consider a model describing the movements of a flock of birds. We take $d = 6$ and $m = 3$. The state is $x=(p,v)\in\mathbb{R}^3\times\mathbb{R}^3$, representing position $p$ and velocity $v$. The dynamics are defined by the coupled stochastic differential equations $dp_t = v_t\,dt$ and $dv_t = \alpha_t\,dt + \sigma\,dW_t$. Note that the velocity dynamics are implicitly coupled with the population distribution through the running cost, since the optimal control depends on the alignment term which involves $\mu$. More generally, the formulation in Section~\ref{sec:MFC} allows the drift $b$ to depend explicitly on $\mu$ through a kernel convolution; the present flocking model illustrates a case where this coupling appears through the cost rather than the drift.
The running cost, promoting velocity alignment, is:
$$f(t,a,p,v,\mu) = \frac{1}{2}\|a\|^2 + \frac{\lambda}{2}\left\| \int_{\mathbb{R}^6}\frac{v-v'}{1+\|p-p'\|^2}\,\mu(dp',dv') \right\|^2\!\!\!.$$
The interaction kernel, which depends on relative distance, is the Cauchy kernel $K(p-p') = \frac{1}{1+\|p-p'\|^2}$. More generally, one may employ the generalized Cauchy kernel, of the form $K(p-p')=\frac{1}{\left(1+\left\Vert p-p'\right\Vert ^{\alpha}\right)^{\beta}}$, where $\alpha\in(0,2]$ and $\beta>0$ \citep{gneiting2004stochastic}.
The terminal cost for this flocking model is
$g(p)=\min(\|p-p_L\|, \|p-p_R\|),$
where $p_L$ and $p_R$ are two target positions. Each agent is penalized by its distance to the nearest target.

\subsection{$N$-Agents Control Problem}\label{sec:n_particle}

MFC problems can naturally be interpreted as limit counterparts of large $N$-agents control problems. Here, we also introduce the corresponding $N$-agents control problems with kernel interactions. Unlike MFC problems, each agent $i\in[N]$ in a $N$-agents problem has their own control $\bm{\alpha}^{i}$, and together form a control profile $\underline{\bm{\alpha}}=(\bm{\alpha}^{1},\ldots,\bm{\alpha}^{N})$ in an appropriate set $\cA^N$. Under a given $\underline{\bm{\alpha}}\in\cA^N$, the dynamics of the controlled $N$-agents is a coupled SDE system:
\begin{equation}\label{eq:4}
dX^{\underline{\bm{\alpha}},i}_t = \tilde{b}\big(t, X^{\underline{\bm{\alpha}},i}_t, \alpha^{i}_t, (K\star\mu^{\underline{\bm{\alpha}}}_t)(X^{\underline{\bm{\alpha}},i}_t)\big)\,dt + \sigma\,dW^i_t,
\end{equation}
for $i\in [N]$, where $\bm{W}^i=(W_t^{i})_{t \in [0,T]}, i\in [N]$, are independent $m$-dimensional Brownian motions and $\mu^{\underline{\bm{\alpha}}}_t=\frac 1 N\sum_{i=1}^N\delta_{X^{\underline{\bm{\alpha}},i}_t}$
is the empirical measure of the state at time $t$. Hence the kernel interaction for agent $i$ at time $t$ reads 
\[
(K\star\mu^{\underline{\bm{\alpha}}}_t)(X^{\underline{\bm{\alpha}},i}_t) = \frac 1 N\sum_{j=1}^N K(X^{\underline{\bm{\alpha}},i}_t-X^{\underline{\bm{\alpha}},j}_t).
\]
The corresponding objective functional is the average over the total $N$ agents:
\begin{equation*}
\begin{split}
J^N(\underline{\bm{\alpha}}) =& 
\frac 1 N\sum_{i=1}^N\mathbb{E}\bigg[\int_0^T \tilde{f}\big(t, X_t^{\underline{\bm{\alpha}},i}, \alpha^{i}_t, (K\star\mu^{\underline{\bm{\alpha}}}_t)(X^{\underline{\bm{\alpha}},i}_t)\big)\,dt\\
& + \tilde{g}\big(X_T^{\underline{\bm{\alpha}},i},(K\star\mu^{\underline{\bm{\alpha}}}_T)(X^{\underline{\bm{\alpha}},i}_T)\big)
\bigg].
\end{split}
\end{equation*}

Hence, the \textbf{$N$-agents control problem} is to seek $\underline{\bm{\alpha}}^\ast\in\cA^N$ such that 
% \begin{equation}\label{eq:6}
$J^N(\underline{\bm{\alpha}}^\ast)=\inf_{\underline{\bm{\alpha}}\in \cA^N} J^N(\underline{\bm{\alpha}}).$
% \end{equation}
% Under suitable assumptions (see, e.g., \citealt{lacker2017limit,djete2022mckean} and~\citealt[Section 6.1.3]{carmona2018probabilistic2} for precise assumptions and rigorous proofs)
% When the drift coefficient $(t,x,\alpha,\mu)\mapsto b(t,x,\alpha,\mu)$ is affine in $x,\alpha,\mu$ and under some convexity and regularity assumptions for the cost functions $f,g$
Under some regularity conditions for the drift coefficient $b$, cost functions $f,g$, and additional convexity assumption for $f,g$ (see the detailed statements in Assumption \textbf{Control of MKV Dynamics} in~\citep[Section 6.4.1]{carmona2018probabilistic})), the optimal control rule for the MFC problem can be shown to serve as an approximate optimal control rule (all agents use the same control rule) for the $N$-agents problem, with $N$ large enough (see Theorem 6.16 in~\citep{carmona2018probabilistic2}). Furthermore, the control of the MFC problem has the advantage of being decentralized in the sense that it depends only on the individual state and not on the states of all agents.  Mathematically, the optimal control of the MFC problem is in feedback form, which is given by a function $a^\ast: [0,T]\times\RR^d\to\mathbb{A}$. This function $a^*$ can be used to construct a decentralized control strategy for $N$-agents problem, which is approximately optimal in the sense that 
$$
    J^N(\bm{\alpha}^{*,1},\ldots,\bm{\alpha}^{*,N})\leq\inf_{\underline{\bm{\alpha}}\in \cA^N} J^N(\underline{\bm{\alpha}})+\epsilon(N),
$$
where for each $i\in[N],
\alpha^{*,i}_t=a^\ast(t,X^i_t), t\in[0,T]$, and $\epsilon(N)\to 0$ as $N\to\infty$.
Note that in general, a feedback control for each agent is a function $\alpha:[0,T]\times(\RR^d)^N\to\mathbb{A}$, which depends on the states of all agents.

\section{Efficient Computation via Random Fourier Features (RFF)}
\label{sec:RFF}

We focus on the interaction term. We will simulate the law of state $X_t$ by using a particle approximation
$\mu_t \approx \frac{1}{N}\sum_{j=1}^N \delta_{X_t^j}$. Note that the goal is to learn the true MFC optimal control rather than solving an $N$-agent control problem.
Computing the convolution term
\begin{equation}\label{eq:convolution}
(K\star\mu_t)(X_t^i) = \frac{1}{N}\sum_{j=1}^N K(X_t^{i}-X_t^{j})
\end{equation}
for every agent $i\in[N]$ requires $\mathcal{O}(N^2)$ operations, which is prohibitive when the number of agents is large.
If $K$ is isotropic, i.e. $K(x-y)=\kappa(\Vert x-y\Vert)$ for some function $\kappa:\mathbb{R}\rightarrow\mathbb{R}$, and positive definite, Bochner’s theorem ensures the existence of a spectral density $p_K(\omega)$ such that
\[
    K(x-y)=K(0)\int_{\mathbb{R}^d} \cos{(\omega\cdot(x-y))}\,p_K(\omega)\,d\omega.
\]
The idea of \textbf{Random Fourier Features} (\textbf{RFF}, \citealt{rahimi2007random}) is to approximate this integral by Monte Carlo sampling:
draw $\omega_m\sim p_K$, $m\in[M]$, and define
\[
    \phi_m(x)=\sqrt{\tfrac{K(0)}{M}}\cos(\omega_m\cdot x), 
    \ \ 
    \phi_{m+M}(x)=\sqrt{\tfrac{K(0)}{M}}\sin(\omega_m\cdot x).
\]
Then the random Fourier feature mapping
\begin{equation}\label{eq:feature_mapping}
    \Phi(x)=(\phi_1(x),\dots,\phi_{2M}(x))^\top
\end{equation}
approximates the kernel $K$:
\(
    K(x-y)\approx \Phi(x)^\top\Phi(y),
\)
so that
\(
    (K\star\mu_t)(X_t^{i})\approx \Phi(X_t^{i})^\top
    \left(
    \frac{1}{N}\sum_{j=1}^N \Phi(X_t^{j})
    \right).
\)
Hence the convolution \eqref{eq:convolution} can be computed for all agents $i\in[N]$ in $\mathcal{O}(NM)$ operations rather than $\mathcal{O}(N^2)$, with approximation error decaying as $\mathcal{O}(M^{-1/2})$.

\section{Applications}
\label{sec:appli}

% \subsection{Gaussian Interactions for Pedestrians}

\noindent\textbf{Gaussian Interactions for Pedestrians. }
For the pedestrian model of Example~1 (Section~\ref{sec:model-crowd}), we use the Gaussian kernel, scaled such that $K(0)=1$,
\(
K(x-y)=\exp\!\left(-\frac{\|x-y\|^2}{2\sigma_K^2}\right).
\)
It is well known that its spectral density is also Gaussian:
\(
p_K(\omega) = \left(\frac{\sigma_K^2}{2\pi}\right)^{d/2}\exp\!\left(-\sigma_{K}^{2}\frac{\|\omega\|^2}{2}\right).
\)
Therefore, we can sample $\omega_m\sim\mathcal{N}(0,{(\sigma_K)}^{-2}I_d)$ for $m\in[M]$, and compute
\(
    (K\star\mu_t)(X_t^{i})\approx 
    \Phi(X_t^{i})^\top
    \Big(\frac{1}{N}\sum_{j=1}^N \Phi(X_t^{j})\Big),
\)
where $\Phi$ is the feature mapping defined in equation~\eqref{eq:feature_mapping}. 
This yields an efficient and scalable approximation of the congestion term 
$\phi((K\star\mu_t)(x))$ in the running cost, enabling simulation or optimization of large-scale pedestrian dynamics in the MFC framework. We omit the computations for the Mat{\'e}rn kernel due to space constraints.

% \subsection{Flocking}

\noindent\textbf{Flocking. }
% \noindent\textbf{Kernel.} 
Here we revisit the flocking model of Example~2 (Section~\ref{sec:model-flocking}). For a single agent with position $p \in \mathbb{R}^{3}$ and velocity $v \in \mathbb{R}^{3}$, consider the alignment vector
\begin{equation}
    A(p,v)
    \;=\;
    \int_{\mathbb{R}^{6}}
        \frac{v - v'}{1+\|p-p'\|^2}\,
        \mu(dp',dv')\in\mathbb{R}^{3}.
\end{equation}
The corresponding running cost is the squared norm $\|A(p,v)\|^2.$
Let us define the isotropic Cauchy kernel
\begin{equation}\label{eq:cauchy_kernel}
    K(p-p') \;=\; \frac{1}{1+\|p-p'\|^2}.
\end{equation}
Let $\otimes$ denote the Hadamard element-wise product. 
Then 
% \begin{align*}
$
    A(p,v)
    % &= v\int_{\mathbb{R}^6} K(p-p')\,\mu(dp',dv') 
    % \\
    % &\quad - \int_{\mathbb{R}^6} v' K(p-p')\,\mu(dp',dv') \notag\\[4pt]
    = v\otimes\mathcal{C}_0(p) - \mathcal{C}_1(p),
$
% \end{align*}
with
$
    \mathcal{C}_\ell(p) := \int_{\mathbb{R}^6} (v')^\ell K(p-p')\,\mu(dp',dv')\in\mathbb{R}^{3},
$
where the power $\ell$ is applied element-wise. Thus
$
    \|A(p,v)\|^2 = \|v\otimes\mathcal{C}_0(p) - \mathcal{C}_1(p)\|^2.
$

% \noindent\textbf{Empirical (particle) approximation.}
Let $\mu$ be approximated by an empirical measure with $N$ particles
$(p_i,v_i)_{i=1}^N$.
Then
\begin{equation*}
    \mathcal{C}_\ell(p) \approx \frac{1}{N}\sum_{j=1}^N (v_j)^\ell\,K(p-p_j)\in\mathbb{R}^{3}.
\end{equation*}
At the particle position $p=p_i$ we have
\begin{align*}
    A_i 
    &:= A(p_i,v_i)\approx
    %  \frac{1}{N}\sum_{j=1}^N  (v_i - v_j)\,K(p_i-p_j)
    % \notag\\[4pt]
    % &\,= 
    v_i  \frac{1}{N} \sum_{j=1}^N  K(p_i-p_j)
    -  \frac{1}{N} \sum_{j=1}^N  v_j K(p_i-p_j).
\end{align*}
Define 
$
    \mathcal{C}_{N,\ell,i} :=  \frac{1}{N}\sum_{j=1}^N (v_j)^\ell K(p_i-p_j)\in\mathbb{R}^{3},
$
so that
$
    A_i \approx v_i \otimes \mathcal{C}_{N,0,i} - \mathcal{C}_{N,1,i}.
$
The empirical running cost is then
$
    \frac{1}{N}\sum_{i=1}^N \|A_i\|^2.
$

% \noindent\textbf{RFF factorization of the kernel.} 
Since the Cauchy kernel \eqref{eq:cauchy_kernel} is isotropic and positive definite, it admits the random Fourier features approximation
$
    K(p-p') \;\approx\; \Phi(p)^\top \Phi(p'),
$
defined in equation~\eqref{eq:feature_mapping}. Its spectral density can be sampled via a Gamma--Gaussian scale mixture: draw $t\sim\mathrm{Gamma}(\beta,1)$ and set $\omega=\sqrt{2t}\,z$ with $z\sim\mathcal{N}(0,I_d)$\citep{langrene2025mixture}. Then
\begin{align*}
\mathcal{C}_{N,\ell,i} &\approx
\frac{1}{N} \sum_{j=1}^N (v_j)^\ell\,\Phi(p_j)^\top \Phi(p_i) = S_\ell\, \Phi(p_i),
\end{align*}
where the feature-domain aggregate $S_\ell$ is defined by
\begin{equation*}
    S_\ell :=  \frac{1}{N} \sum_{j=1}^N (v_j)^{\ell}\,\Phi(p_j)^\top \in \mathbb{R}^{3\times 2M}.
\end{equation*}
The empirical alignment vector can be approximated by
\begin{equation*}
    A_i \approx v_i\otimes S_0 \Phi(p_i)
    - S_1 \Phi(p_i) = ((v_i\mathbf{1}_{2M}^\top)\otimes S_0-S_1) \Phi(p_i),
\end{equation*}
and the empirical running cost can be approximated by:
\begin{equation*}
    \frac{1}{N}\sum_{i=1}^N \|A_i\|^2 \approx \frac{1}{N}\sum_{i=1}^N \|((v_i\mathbf{1}_{2M}^\top)\otimes S_0-S_1) \Phi(p_i)\|^2,
\end{equation*}
where we stress that $S_0$ and $S_1$ are independent of $i$.

\section{Numerical results}
\label{sec:numerics}

% \subsection{Pseudo-code}

\textbf{Algorithm.} We use a direct policy search method: We take a feedback control which is a neural network and train it to minimize the total cost; this idea was proposed in~\citep[Algorithm 1]{carmona2022convergence} using an empirical distribution. Algorithm~\ref{alg:mfc-generic} summarizes the approach. For simplicity, $\tilde{g}$ is independent of $\mu_T$ but kernel-based interactions could be included as well. Next, we solve numerical examples and compare the performance of the ``{\bf naive}'' method and our proposed method based on {\bf RFF}.

\begin{algorithm2e}[t]
\footnotesize
\DontPrintSemicolon
Initialize $\theta$ for controller $\alpha_\theta$\;
\For{$k=0$ \KwTo $K$ \textnormal{(SGD iterations)}}{
    Sample $(\Delta W_n^i)_{i=1}^N$ and $(X_0^i)_{i=1}^N$\;
    \For{$n=0$ \KwTo $T/\Delta t -1$ \textnormal{(simulate trajectory)}}{
        $t \gets n\Delta t$; $\alpha_n^i \gets a_\theta(t,X_n^i)$, $i \in [N]$\;
        \eIf{\textnormal{naive}}{
            $\kappa_n^j \gets \frac1N \sum_{i=1}^N K(X_n^j - X_n^i)$, $j \in [N]$\;
        }{
            \textnormal{RFF}: $\phi_n^i \gets \phi(X_n^i)$ $i \in [M]$, $\overline{\phi} \gets \tfrac1N\sum_i \phi_n^i$\;
            $\kappa_n^j \gets \langle \phi_n^j ,\, \overline{\phi} \rangle$, $j \in [N]$\;
        }
        $X_{n+1}^i \gets X_n^i + \tilde b(t,X_n^i,\alpha_n^i,\kappa_n^i)\,\Delta t + \sigma\,\Delta W_n^i$, $i \in [N]$\;
        $c_n^i \gets \tilde f(t,X_n^i,\alpha_n^i,\kappa_n^i)$, $i \in [N]$\;
    }
    $C(\theta) \gets \frac1N \sum_i \left(\sum_n c_n^i\Delta t + \tilde g(X_T^i)\right)$\;
    Update $\theta$ via gradient step on $C(\theta)$\;
}
\caption{SGD for MFC with Kernel Interactions}
\label{alg:mfc-generic}
\end{algorithm2e}

\textbf{Experimental setup.} The feedback control is a fully connected neural network with $\tanh$ activations: two hidden layers with 64 units for crowd motion, three hidden layers with 128 units for flocking. The network takes $(t,x)$ as input and outputs the control action. Training uses Adam with 1000 SGD iterations, learning rate $3\times10^{-4}$, and $N=200$ particles for crowd motion; 5000 iterations, learning rate $5\times10^{-5}$, and $N=300$ particles for flocking. We use $M=512$ random Fourier features throughout. After training, we evaluate control quality by simulating with the exact kernel and $N=1000$ particles (the \emph{evaluation cost}). In the figures below, left panels show the evaluation cost during training and right panels show training time as a function of $N$ (averaged over 10 runs; dashed lines: log-log fit $T\propto N^p$; shaded areas: $\pm1$ standard deviation). In all cases the RFF method preserves control quality while confirming the expected complexity scaling. Experiments ran on an NVIDIA RTX 4090 D GPU (Intel Xeon Gold 6530 CPU, PyTorch 2.0.1); timing used CPU-only nodes.

% \subsection{Crowd motion}

\textbf{Crowd motion.} 
We solve numerically Example~1 (Section~\ref{sec:model-crowd}) with two types of kernels: Gaussian and generalized Mat{\' e}rn. In both cases, the other model parameters are the same and are chosen as: $T=2.0, \sigma=0.1, \lambda=5.0$. For the time discretization, we use $\Delta t=0.05$.

\begin{figure}[htbp]
    \centering
    \begin{subfigure}{0.2\textwidth}
        \includegraphics[width=\textwidth]{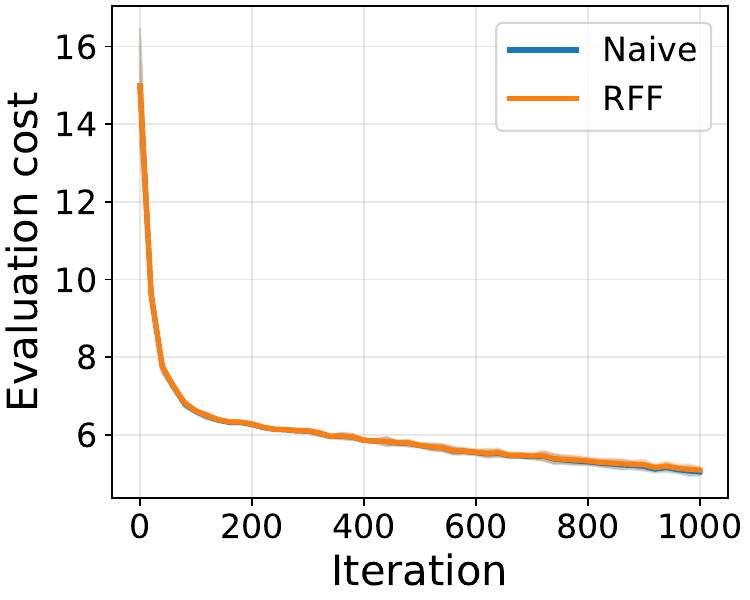}
    \end{subfigure}
    \begin{subfigure}{0.2\textwidth}
        \includegraphics[width=\textwidth]{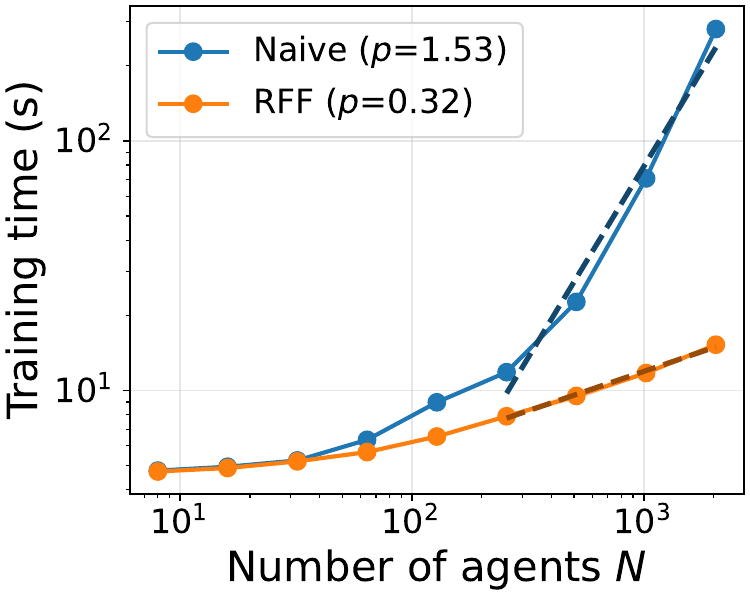}
    \end{subfigure}
    \caption{Ex. 1, Gaussian kernel: cost and training time.}
    \label{fig:crowd-gaussian-cost-compute}
\end{figure}

\textit{Gaussian kernel.} We use a Gaussian kernel with $\sigma_k=0.3$. Fig.~\ref{fig:crowd-gaussian-cost-compute} shows the cost and timing results; the training time scaling is roughly consistent with our theoretical estimates. We omit the density plots due to space constraints.
% Fig.~\ref{fig:crowd-gaussian-density} shows the density evolution under the learned control; the red marker indicates the target.

% \begin{figure}[htbp]
%     \centering
%     \begin{subfigure}{0.2\textwidth}
%         \includegraphics[width=\textwidth]{FIGURES/CROWD-GAUSSIAN/crowd-GaussianKTrue-1run-density-3d-t0.pdf}
%         % \includegraphics[width=\textwidth]{FIGURES/CROWD-GAUSSIAN/crowd-GaussianKTrue-1run-density-3d-t0-uniform.pdf}
%     \end{subfigure}
%     \begin{subfigure}{0.2\textwidth}
%         \includegraphics[width=\textwidth]{FIGURES/CROWD-GAUSSIAN/crowd-GaussianKTrue-1run-density-3d-t1.pdf}
%         % \includegraphics[width=\textwidth]{FIGURES/CROWD-GAUSSIAN/crowd-GaussianKTrue-1run-density-3d-t1-uniform.pdf}
%     \end{subfigure}\\
%     \begin{subfigure}{0.2\textwidth}
%         \includegraphics[width=\textwidth]{FIGURES/CROWD-GAUSSIAN/crowd-GaussianKTrue-1run-density-3d-t2.pdf}
%         % \includegraphics[width=\textwidth]{FIGURES/CROWD-GAUSSIAN/crowd-GaussianKTrue-1run-density-3d-t2-uniform.pdf}
%     \end{subfigure}
%     \begin{subfigure}{0.2\textwidth}
%         \includegraphics[width=\textwidth]{FIGURES/CROWD-GAUSSIAN/crowd-GaussianKTrue-1run-density-3d-t4.pdf}
%         % \includegraphics[width=\textwidth]{FIGURES/CROWD-GAUSSIAN/crowd-GaussianKTrue-1run-density-3d-t4-uniform.pdf}
%     \end{subfigure}
%     \caption{Ex. 1, Gaussian kernel: density evolution.}
%     \label{fig:crowd-gaussian-density}
% \end{figure}

\begin{figure}[htbp]
    \centering
    \begin{subfigure}{0.2\textwidth}
        \includegraphics[width=\textwidth]{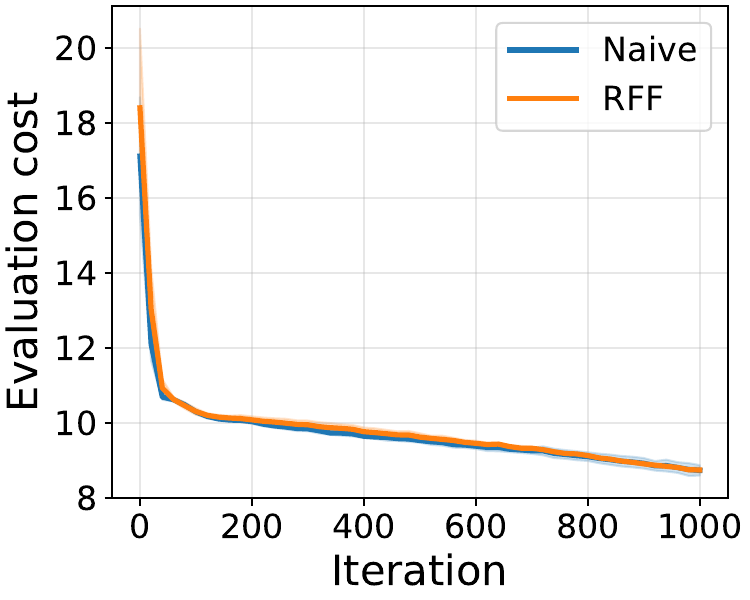}
    \end{subfigure}
    \begin{subfigure}{0.2\textwidth}
        \includegraphics[width=\textwidth]{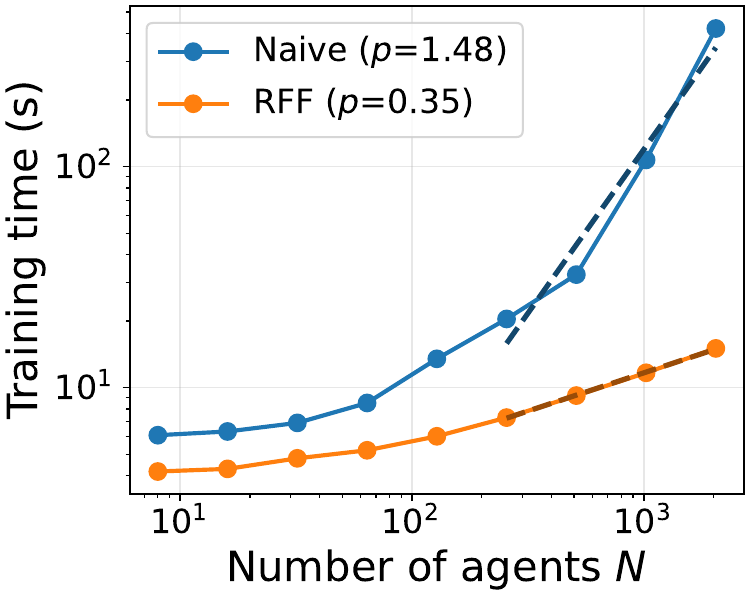}
    \end{subfigure}
    \caption{Ex. 1, Mat{\' e}rn kernel: cost and training time.}
    \label{fig:crowd-matern-cost-compute}
\end{figure}

\begin{figure}[htbp]
    \centering
    \begin{subfigure}{0.2\textwidth}
        \includegraphics[width=\textwidth]{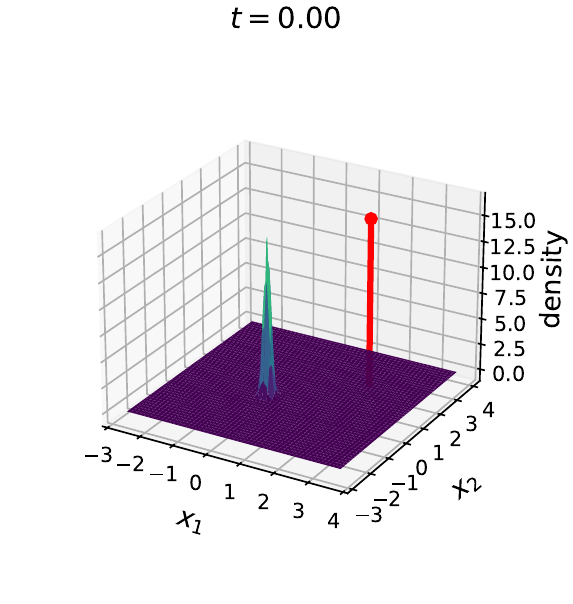}
    \end{subfigure}
    \begin{subfigure}{0.2\textwidth}
        \includegraphics[width=\textwidth]{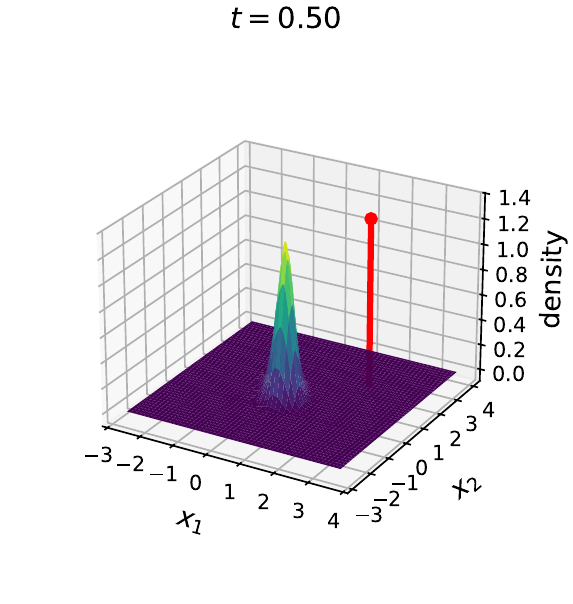}
    \end{subfigure}\\
    \begin{subfigure}{0.2\textwidth}
        \includegraphics[width=\textwidth]{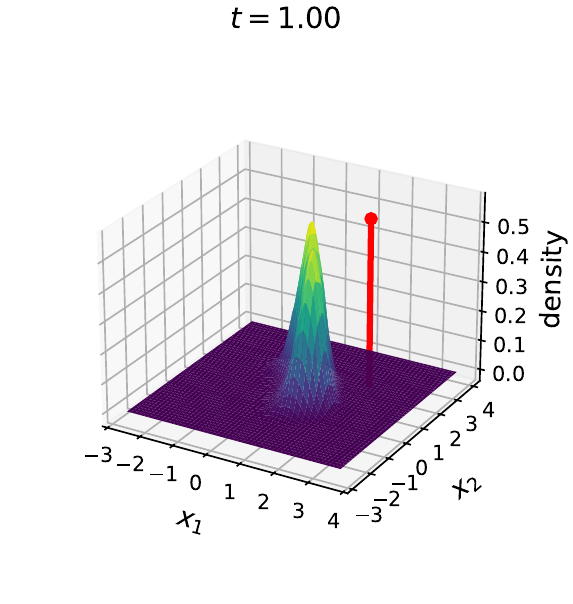}
    \end{subfigure}
    \begin{subfigure}{0.2\textwidth}
        \includegraphics[width=\textwidth]{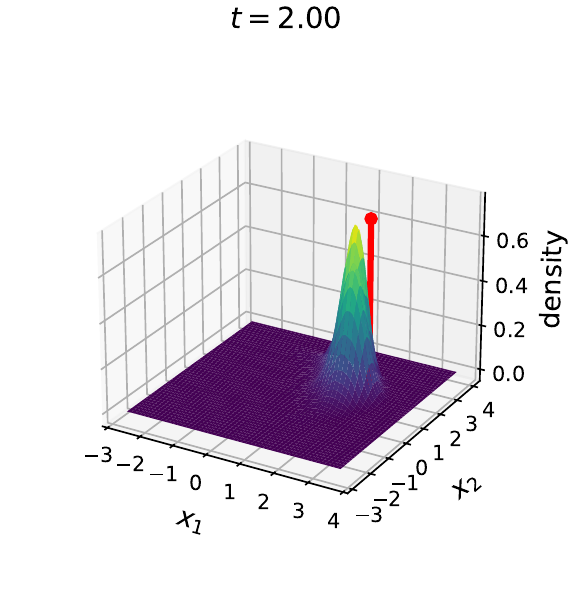}
    \end{subfigure}
    \caption{Ex. 1, Mat{\' e}rn kernel: density evolution.}
    \label{fig:crowd-matern-density}
\end{figure}

\begin{figure}[htbp]
    \centering
    \begin{subfigure}{0.2\textwidth}
        \includegraphics[width=\textwidth]{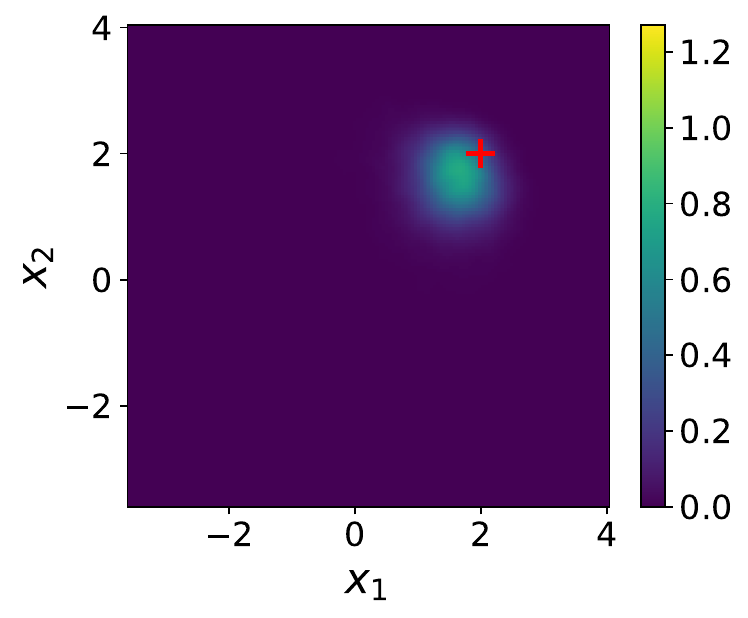}
        \caption{Gaussian}
    \end{subfigure}
    \begin{subfigure}{0.2\textwidth}
        \includegraphics[width=\textwidth]{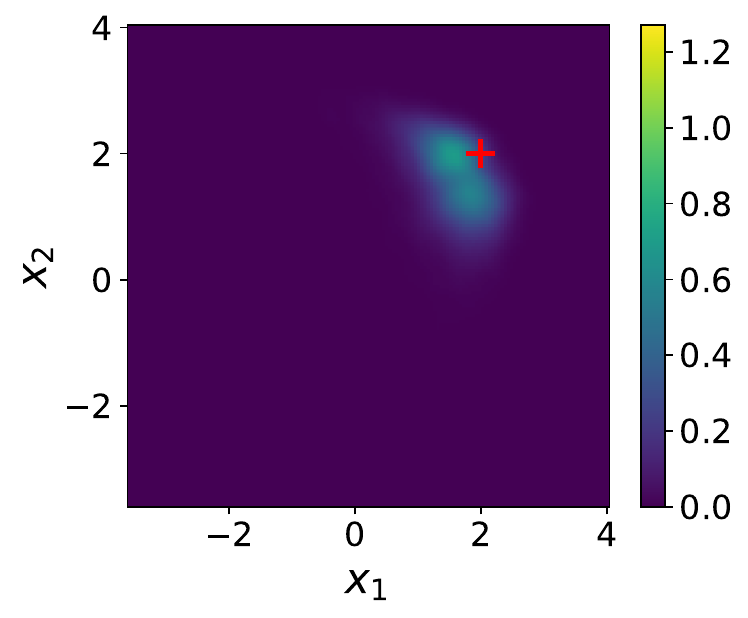}
        \caption{$\beta_1\!=\!1.9,\,\beta_2\!=\!1.5$}
    \end{subfigure}
    \begin{subfigure}{0.2\textwidth}
        \includegraphics[width=\textwidth]{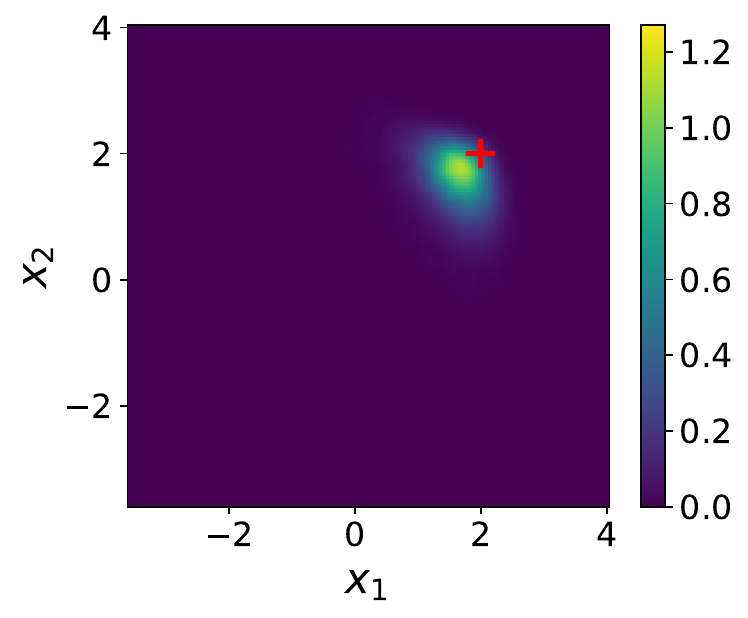}
        \caption{$\beta_1\!=\!1.0,\,\beta_2\!=\!1.5$}
    \end{subfigure}
    \begin{subfigure}{0.2\textwidth}
        \includegraphics[width=\textwidth]{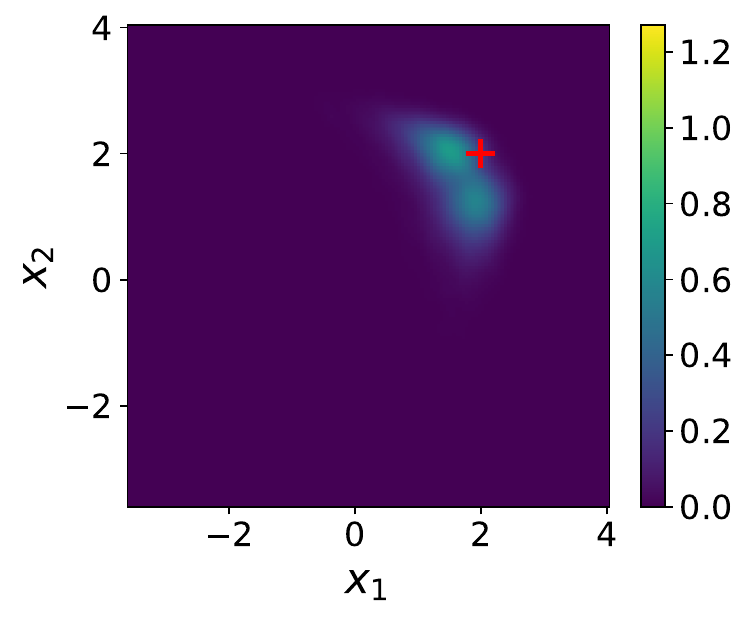}
        \caption{$\beta_1\!=\!1.9,\,\beta_2\!=\!5.0$}
    \end{subfigure}
    \caption{Ex.~1: terminal density ($t=T$) for different kernels. (a)~Gaussian ($\sigma_k=0.3$); (b)--(d)~generalized Mat{\'e}rn with varying $\beta_1$ (tail decay) and $\beta_2$ (smoothness). All panels share the same color scale.}
    \label{fig:crowd-gaussian-density-heatmap-cmp}
\end{figure}

\textit{Generalized Mat{\'e}rn kernel.} We use $\beta_1=1.9$, $\beta_2=1.5$. Fig.~\ref{fig:crowd-matern-cost-compute} shows the cost and timing results. Fig.~\ref{fig:crowd-matern-density} shows the density evolution. Fig.~\ref{fig:crowd-gaussian-density-heatmap-cmp} compares the terminal density across four kernel configurations: Gaussian and three generalized Mat{\'e}rn kernels with different $\beta_1$ (tail decay) and $\beta_2$ (smoothness). Decreasing $\beta_1$ from $1.9$ to $1.0$ yields heavier tails, while increasing $\beta_2$ from $1.5$ to $5.0$ produces a smoother kernel that approaches the Gaussian limit.

\begin{figure}[htbp]
    \centering
    \begin{subfigure}{0.2\textwidth}
         \includegraphics[width=\textwidth]{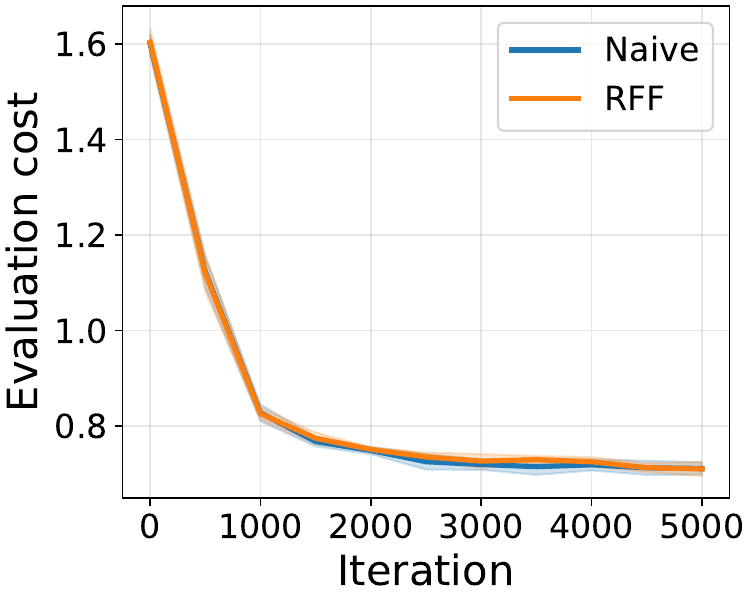}
    \end{subfigure}%
    \begin{subfigure}{0.2\textwidth}
         \includegraphics[width=\textwidth]{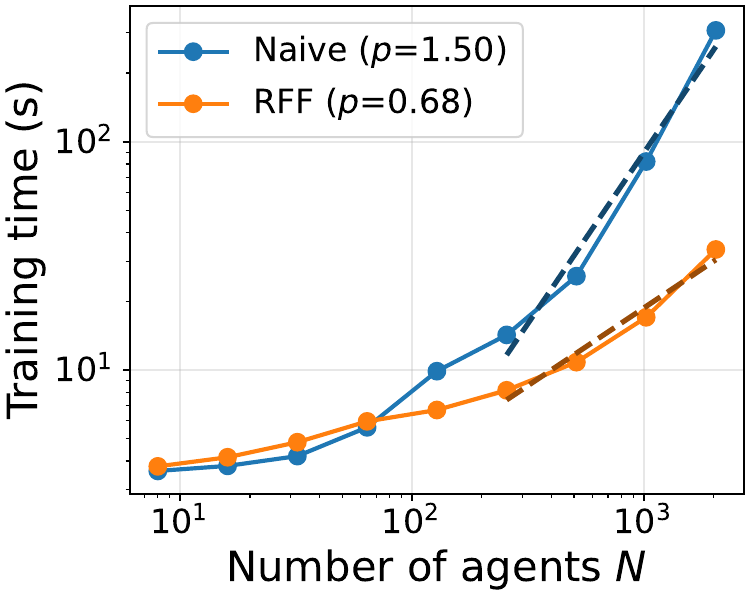}
    \end{subfigure}
    \caption{Ex. 2, Cauchy kernel: cost and training time.}
    \label{fig:flock-cauchy-cost-compute}
\end{figure}

\begin{figure}[htbp]
    \centering
    \begin{subfigure}{0.2\textwidth}
        \includegraphics[width=\textwidth]{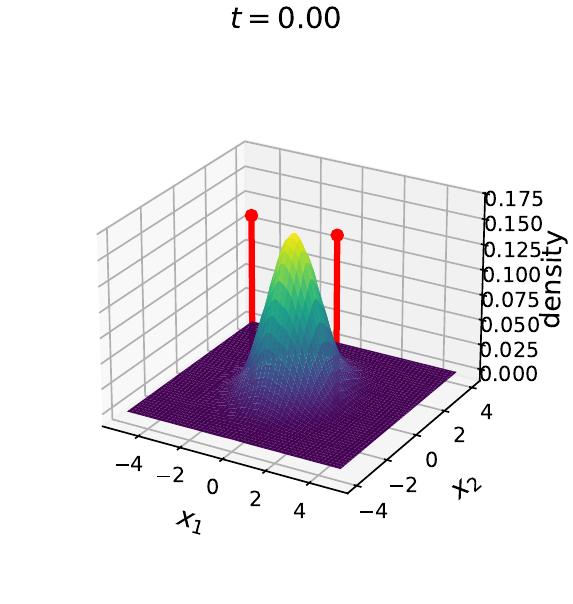}
    \end{subfigure}
    \begin{subfigure}{0.2\textwidth}
        \includegraphics[width=\textwidth]{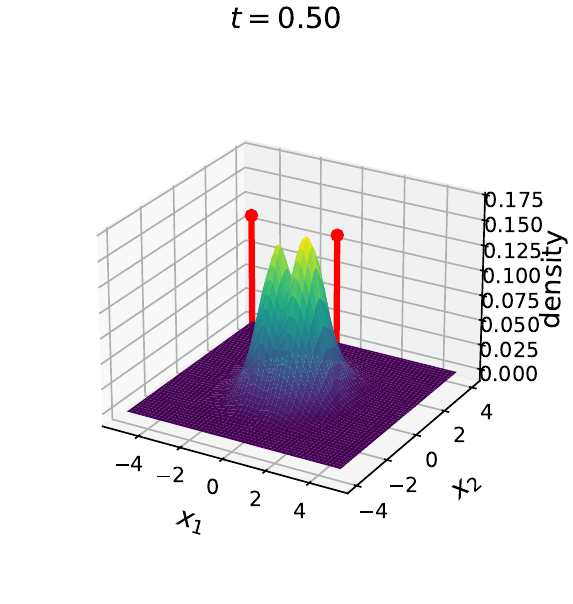}
    \end{subfigure}\\
    \begin{subfigure}{0.2\textwidth}
        \includegraphics[width=\textwidth]{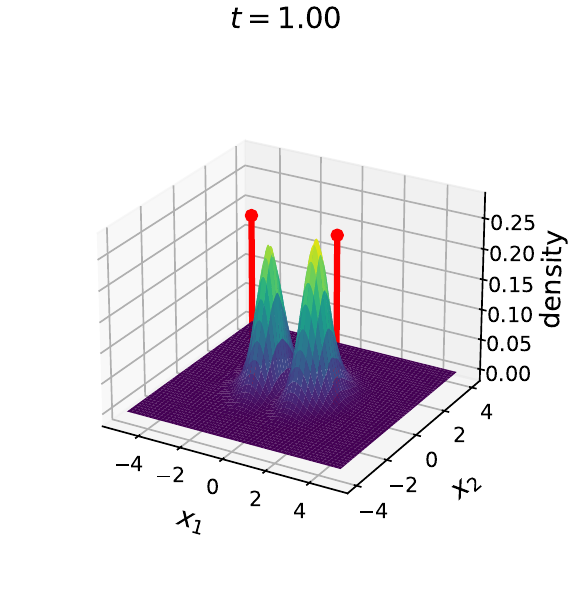}
    \end{subfigure}
    \begin{subfigure}{0.2\textwidth}
        \includegraphics[width=\textwidth]{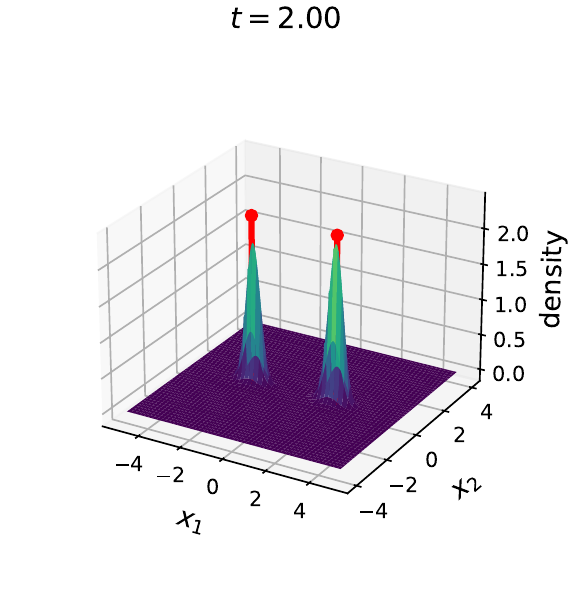}
    \end{subfigure}
    \caption{Ex. 2: position's density evolution.}
    \label{fig:flock-cauchy-density-pos}
\end{figure}

\begin{figure}[htbp]
    \centering
    \begin{subfigure}{0.2\textwidth}
        \includegraphics[width=\textwidth]{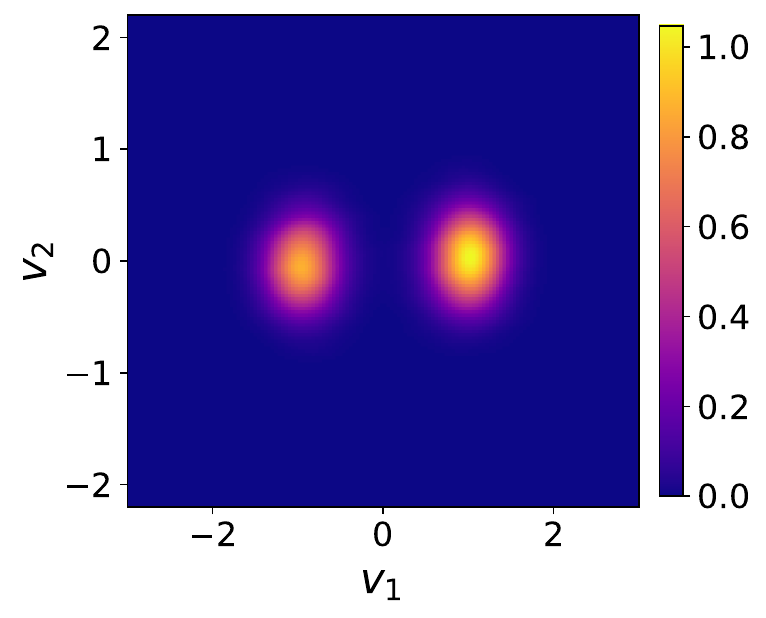}
    \end{subfigure}
    \begin{subfigure}{0.2\textwidth}
        \includegraphics[width=\textwidth]{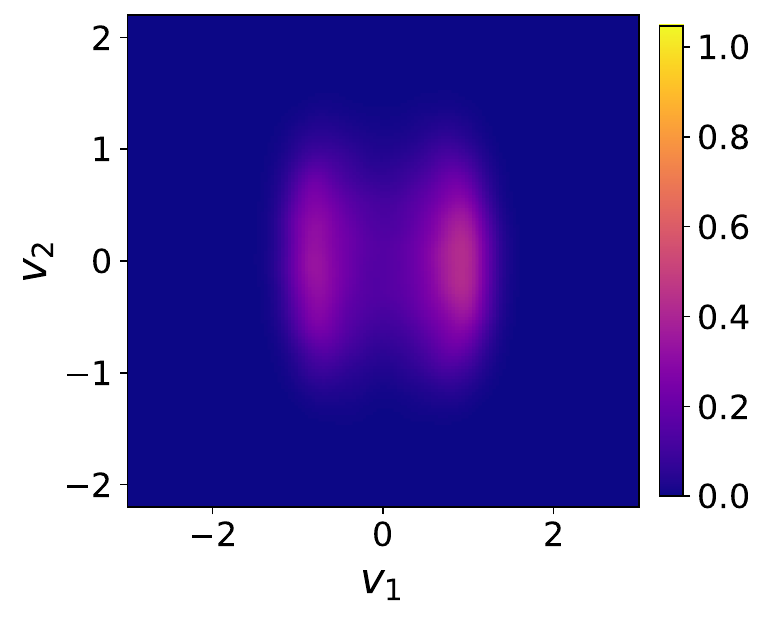}
    \end{subfigure}
    \caption{Ex. 2: velocity density at $t=T/2$ with $\beta=1$ (left) and $\beta=10$ (right).}
    \label{fig:flock-cauchy-density-cmp-beta}
\end{figure}

% \subsection{Flocking}
\textbf{Flocking.}
We now turn to Example~2 (Section~\ref{sec:model-flocking}) about flocking, using the generalized Cauchy kernel $(1+\|p-p'\|^2)^{-\beta}$. We solve the problem in dimension $2+2=4$ to visualize the density. The agents start around $(0,0)$ and are attracted to two targets at $(\pm2,0)$, with terminal cost $g(p)=\min(\|p-p_L\|,\|p-p_R\|)$. We use $\beta=10$ and $T=2$, so agents focus on close neighbors. Fig.~\ref{fig:flock-cauchy-cost-compute} shows the cost and timing results. Fig.~\ref{fig:flock-cauchy-density-pos} shows the position density evolution. With $\beta=1$ (same $T=2$), agents align velocity across the whole population and move toward a single target, demonstrating that the behavioral difference is driven by $\beta$ (interaction range) rather than the time horizon. Fig.~\ref{fig:flock-cauchy-density-cmp-beta} shows the velocity density at $t=T/2$: with $\beta=1$ (long-range), two sharply separated velocity modes emerge, whereas with $\beta=10$ (short-range), the velocity distribution remains more diffuse.

\section{Conclusion}
\label{sec:ccl}

We introduced a scalable method for solving MFC problems with kernel interactions by combining stochastic gradient descent with random Fourier features. The main contribution is the replacement of the empirical kernel computation, which scales quadratically in the number of agents, by an approximate representation with linear complexity, including for non-Gaussian kernels. This yields a practical approach for training feedback controls in settings where large interacting populations must be simulated repeatedly.
The numerical experiments illustrate the main computational advantage of the method. Across all examples, including crowd motion and flocking, the RFF-based approximation preserves the quality of the learned control while significantly reducing runtime as the population size increases.
Several directions merit further investigation, such as a more detailed analysis of approximation errors, a systematic study of the effect of $M$ on accuracy and runtime across kernel families, and extending the method to problems with heterogeneous agents, additional interaction models, and constrained dynamics.

\bibliography{ifacconf}

\end{document}